# A class of non-linear fractional-order system stabilisation via fixed-order dynamic output feedback controller

Elyar Zavary
*Advanced Control Systems Laboratory, School of Electrical and Computer Engineering*
*Tarbiat Modares University*
Tehran, Iran
elyar.zavary@modares.ac.ir

Mahdi Sojoodi*
*Advanced Control Systems Laboratory, School of Electrical and Computer Engineering*
*Tarbiat Modares University*
Tehran, Iran
sojoodi@modares.ac.ir

***Abstract*—This paper investigates the robust stabilisation of a class of fractional-order non-linear systems via fixed-order dynamic output feedback controller in terms of linear matrix inequalities (LMIs). The systematic stabilisation algorithm design for low-order controller based on direct Lyapunov approach is proposed. In the presented algorithm the conditions containing the bilinear variables are decoupled into separate conditions without imposing equality constraints or considering an iterative search of the controller parameters. There is no any limiting constraint on the state space matrices and also we assumed the most complete output feedback controller. Simulations results are given to approve the effectiveness and the straightforwardness of the proposed design.**

## I. Introduction

In recent decades, study of fractional-order systems has been expanded significantly. For example hereditary and long memory attribute of systems, such as viscoelastic polymers [1], biomedical applications [2], semi-infinite transmission lines with losses [3], dielectric polarization [4], have been described with fractional-order operators. More, stability analysis of fractional-order systems have been attracted considerable interests, where several literatures addressing this topic have been released [5–8] and subsequently fractional-order controllers design and implementation in system control field have become commonplace [9–11]. A fractional-order PI/sup /spl lambda//D/sup /spl mu/ controller was proposed in [9]. Stability and stabilisation of fractional-order interval systems are studied in [11].

Design of fractional-order $PI^{\lambda}D^{\mu}$ controllers with an improved differential evolution is proposed in [12]. Using the Lyapunov function method,[13] investigates the design of state feedback stabilization controllers for fractional-order nonlinear systems in triangular form. Necessary and sufficient stability conditions of fractional-order interval linear systems are stablished in [14]. In [14] the necessary and sufficient stability conditions of fractional-order systems are directly extended to the robust stability condition of fractional-order interval polynomial systems. Estimation of the system states and observer-based stabilisation were investigated in [15,16]. In [17] using continuous frequency distribution, the stability conditions of a class of Lipschitz nonlinear fractional-order systems based on indirect approach to Lyapunov stability are derived.

The definition of Mittag-Leffler stability definition was proposed in [18], and also fractional Lyapunov direct method was introduced. By using of Mittag-Leffler function, Laplace transform, and the generalized Gronwall inequality, a new sufficient condition ensuring local asymptotic stability and stabilization of a class of fractional-order nonlinear systems with fractional-order $1 < \alpha < 2$ is proposed in [19]. Stability analysis of fractional-order systems is studied in [20], in [20] an extension of Lyapunov direct method for fractional-order systems is proposed. Moreover the studies of Li, Wang and Lu [10] focuses on the observer-based stability problem of a class of non-linear fractional-order uncertain systems with admissible time-variant uncertainty. The proposed method therein is used for stabilisation of a class of nonlinear fractional-order system by assuming that input matrix of the system is of full row rank.

Note that most of the mentioned works focuses on stability study of linear fractional-order systems in which state feedback control law is the most existed control law. Existence of some technical and economic limitations makes it difficult to obtain the system states in practical applications. Output feedback controller eliminates mentioned problems of control and besides that among output feedback controllers, dynamic ones have more degrees of freedom in controller designing procedure and subsequently satisfying control objectives compared with static ones [21]. Due to achieving control objectives most of dynamic controller design methods lead to high order controllers, where High order controllers are not preferable along of costly implementation and maintenance, high fragility, and potential numerical errors [22]. Given that closed-loop performance could not be guaranteed through order reduction methods, it is worthwhile to have a solution to design a controller with low and fixed-order which can be as small as possible to satisfy control objectives [23-24].

The dynamic output feedback controller is a powerful method for controlling the strict feedback nonlinear systems. To the best of our knowledge, there is few results on designing dynamic output feedback controller for the stability of nonlinear fractional-order systems in the literature, this motivated us for the study of this paper. This paper investigates the fixed predetermined order dynamic output feedback controller for the robust stabilisation of fractional-

*Corresponding Author, Address: P.O. Box 14115-…, Tehran, Iran, E-mail address: *sojoodi@modares.ac.ir*, web page: *http://www.modares.ac.ir/~sojoodi* , Tel./Fax: +98 21 8288-3902.

order nonlinear systems with Lipschitz nonlinearities in the states and inputs. It should be mentioned that nonlinear structure used in this paper is commonplace in many real systems [11,25]. No limiting constraints on the state space matrix, assumed in [10], are considered and also the most complete model of dynamic output feedback controller is taken. Notwithstanding this, results are given in terms of linear matrix inequalities (LMIs) where the design parameters can be easily obtained by accessing the feasibility of LMI constraints through optimisation parsers and solvers.

The rest of this paper organised as follows: in Section 2, some preliminaries and problem formulation are presented. The proposed fixed-order dynamic output feedback controller with the design algorithm of controller for robust stabilisation of nonlinear fractional-order systems are derived in Section 3. Some numerical examples are provided in Section 4 to illustrate the effectiveness of proposed method. Eventually section 5 draws the conclusion.

## II. Problem Formulation and Preliminaries

Some mathematical notations that are used throughout this paper, are defined here. $A \otimes B$ represents the Kronecker product of matrices $A$ and $B$. The transpose of $M$ is denoted by $M^T$ and $sym(M)$ stands for $M + M^T$. The notation $\star$ denotes symmetric component in matrix. Psuedo inverse of a given non-square matrix $A_{n\times m}$ is shown by $A^{\uparrow}$.

Fractional-order nonlinear system with the following dynamic is considered

$$\begin{aligned} D^q x(t) &= \tilde{\mathcal{A}} x(t) + \tilde{\mathcal{B}} u(t) + \phi(x(t), u(t)) \\ y(t) &= \mathcal{C} x(t) \end{aligned} \tag{1}$$

where

$$\begin{aligned} \tilde{\mathcal{A}} &= \mathcal{A} + \Delta\mathcal{A} \\ \tilde{\mathcal{B}} &= \mathcal{B} + \Delta\mathcal{B} \end{aligned} \tag{2}$$

with initial condition

$$x(0) = x_0 \tag{3}$$

where $x(t) \in \Re^n, u(t) \in \Re^m, y(t) \in \Re^p$ are pseudo state, input, measured output, respectively. $\mathcal{A} \in \Re^{n\times n}, \mathcal{B} \in \Re^{n\times m}, \mathcal{C} \in \Re^{p\times n}$ are known constant matrices, and $\phi(\cdot) : [\Re^n \times \Re^m] \to \Re^n$, is nonlinear function. $\Delta\mathcal{A} \in \Re^{n\times n}$ and $\Delta\mathcal{B} \in \Re^{n\times m}$ are time-invariant matrices, with parametric uncertainty. $q$ is the fractional derivative order, there are several definitions for fractional-order derivative, among them Grünwald-Letnikov, Riemann-Liouville and Caputo are most commonly referred. Since the initial condition of Caputo definition is similar to integer order ones, as a physical aspect, Caputo definition is used in this with the following definition

$${}_a^C D_t^q = \frac{1}{\Gamma(\bar{n}-\alpha)} \int_a^t (t-\tau)^{\bar{n}-a-1} \left(\frac{d}{d\tau}\right)^{\bar{n}} f(\tau) d\tau,$$

where $\Gamma(\cdot)$ is Gamma function defined by $\Gamma(\epsilon) = \int_0^\infty e^{-t} t^{\epsilon-1} dt$ and $\bar{n}$ is the smallest integer that is equal or greater than $q$.

**Lemma 1** [10] Let $f: \Re_\epsilon \to \Re^n$ be piecewise continuous respect to $t$, where $\Re_\epsilon = \{(t, x): 0 \le t \le a \text{ and } \|x - x_0\| \le b\}$, $f = [f_1, \dots, f_n]^T$, $x \in \Re^n$ and $\|f(t, x)\| \le M$ on $\Re_\epsilon$. Then, there exists at least one solution for the system of fractional differential equations given by

$$D^q x(t) = f(t, x(t)) \tag{4}$$

with the initial condition

$$x(0) = x_0 \tag{5}$$

on $0 \le t \le \beta$ where $\beta = min\left(a, \left[(b/M)\Gamma(q+1)^{1/q}\right]\right), 0 < q < 1$.

**Lemma 2** [10] Consider initial fractional problem (4) and (5) with $0 < q < 1$ and assume that Lemma 1 conditions hold. Let

$$\begin{aligned} & g(v, x_*(v)) = \\ & f\left(t - \left(t^q - v\Gamma(q+1)\right)^{1/q}, x\left(t - v\Gamma(q+1)\right)^{1/q}\right) \end{aligned}$$

then $x(\mathrm{t})$, is given by

$$x(t) = x_*\left(t^q / \Gamma(q+1)\right),$$

where $x_*(\upsilon)$ can be obtained by solving the following integer order differential equation

$$\begin{aligned} \frac{dx_*(v)}{dv} &= g(v, x_*(v)) \\ x(0) &= x_0. \end{aligned} \tag{6}$$

System matrices $\mathcal{A}$, $\mathcal{B}$, $\mathcal{C}$, nonlinear function $\phi(\cdot)$ and uncertainty matrices $\Delta\mathcal{A}$ and $\Delta\mathcal{B}$ are assumed to satisfy the following assumptions.

**Assumption 1.** The pairs of $(\mathcal{A}, \mathcal{B})$ and $(\mathcal{A}, \mathcal{C})$ are controllable and observable, respectively.

**Assumption 2.** $\Delta\mathcal{A}$ and $\Delta\mathcal{B}$ are time-invariant matrix of the following form:

$$[\Delta\mathcal{A} \quad \Delta\mathcal{B}] = \mathcal{M}\Delta(\sigma)[\mathcal{N}_1 \quad \mathcal{N}_2] \tag{7}$$

$$\Delta(\sigma) = \mathcal{Z}(\sigma)[I + \mathcal{J}\mathcal{Z}(\sigma)]^{-1} \tag{8}$$

$$Sym\{\mathcal{J}\} > 0, \tag{9}$$

where $\mathcal{M} \in \Re^{n\times m_0}$, $\mathcal{N}_1 \in \Re^{m_0\times n}$, $\mathcal{N}_2 \in \Re^{m_0\times m}$ and $\mathcal{J} \in \Re^{m_0\times m_0}$ are real known matrices. The uncertain matrix $\mathcal{Z}(\sigma) \in \Re^{m_0\times m_0}$ satisfies

$$Sym\{\mathcal{Z}(\sigma)\} \ge 0,, \tag{10}$$

where $\sigma \in \Omega$, with $\Omega$ being a compact set.

**Remark 1.** Condition *(9)* guarantees that $I + \mathcal{J}\mathcal{Z}(\sigma)$ is invertible for all $\mathcal{Z}(\sigma)$ satisfying *(10)*. Therefore $\Delta(\sigma)$ in *(7)* is well defined ([8]).

**Assumption 3.** Nonlinear function $\phi(x(\mathrm{t}), u(\mathrm{t}))$ is Lipschitz on $x(\mathrm{t})$ with Lipschitz constant $\xi$

$$\begin{aligned} & \|\phi(x_1(t), u_1(t)) - \phi(x_2(t), u_2(t))\| \\ & \qquad < \xi \|x_1(t) - x_2(t)\| \end{aligned} \tag{11}$$

for all $x_1(t), x_2(t) \in \Re^n$ and

$$\phi(0, 0) = 0. \tag{12}$$

**Lemma 3** [6] Let $\mathcal{A} \in \Re^{n\times n}, 0 < q < 1$ and $\theta = (1 - q\pi)/2$. The fractional-order system $D^q x(t) = \mathcal{A}x(t)$ is asymptotically stable if and only if there exist a positive definite Hermitian matrices $X = X^* > 0, X \in \mathfrak{C}^{n\times n}$ such that

$$(rX + \bar{r}\bar{X})^T A^T + A(rX + \bar{r}\bar{X}) < 0, \tag{13}$$

Where $r = e^{\theta i}$.

Notations: In this paper $A \otimes B$ denotes the kronecker product of matrices A and B, and the symmetric of matrix $M$ will be shown by $sym(.)$, which is defined by $sym(M) = M^T + M$, and also $\uparrow$ is the symbol of pseudo inverse of matrix.

**Lemma 4** [7] Let
$\Omega = \{\Delta \in \Re^{m_0\times m_0} | \Delta \text{ is subjected to } (7) - (9)\}$. Then
$\Omega = \{\Delta \in \Re^{m_0\times m_0} | \det(I - \Delta\mathcal{J}) \ne 0 \text{ and } \Delta Sym\{\mathcal{J}\}\Delta^T \le Sym\{\Delta\}$.

## III. Main Result

In this work we will study the stability and asymptotically stabilisation of FOMASs composed of *(1)*, with fixed-order dynamic output feedback controller.

In order to achieve the objectives on system (1), we use the following non-fragile control protocol

$$\begin{aligned} &D^q x_c(t) = \mathcal{A}_c x_c(t) + \mathcal{B}_c \mathcal{Y}(t), \\ &u = \mathcal{C}_c x_c + \mathcal{D}_c \mathcal{Y}(t), \\ &x_c(0) = x_{c0} \end{aligned} \tag{14}$$

where $x_c \in \Re^{n_c}$ is controller pseudo state in which $n_c$ is the controller order and $\mathcal{A}_c, \mathcal{B}_c, \mathcal{C}_c$ and $\mathcal{D}_c$ are controller matrices to be designed.

By implementing the controller (14) on the system (1), the closed-loop system is achieved as follows

$$\begin{aligned} &D^q \mathcal{X}(t) = \Phi(\mathcal{X}(t), t) = \mathcal{A}_{cl,\Delta}\mathcal{X}(t) + \begin{bmatrix} \phi(x(t), u(t)) \\ 0 \end{bmatrix}, \\ &\mathcal{X}(0) = \mathcal{X}_0 = [x_0^T \quad x_{c0}^T]^T \end{aligned} \tag{15}$$

where

$$\begin{aligned} &\mathcal{X}(t) = [x^T(t) \quad x_c^T(t)]^T, \ \mathcal{A}_{cl,\Delta} = \mathcal{A}_\psi + \mathcal{A}_\Delta \\ &\mathcal{A}_\psi = \begin{bmatrix} \mathcal{A} + \mathcal{B}\mathcal{D}_c\mathcal{C} & \mathcal{B}\mathcal{C}_c \\ \mathcal{B}_c\mathcal{C} & \mathcal{A}_c \end{bmatrix}, \mathcal{A}_\Delta = \widetilde{\mathcal{M}}\Delta\widetilde{\mathcal{N}} \\ &\widetilde{\mathcal{M}} = [\mathcal{M}^T \quad 0]^T, \widetilde{\mathcal{N}} = [\mathcal{N}_1 + \mathcal{N}_2\mathcal{D}_c\mathcal{C} \quad \mathcal{N}_2\mathcal{C}_c], \end{aligned} \tag{16}$$

**Theorem 1.** Consider the nonlinear fractional-order system (1) with output dynamic controller (14) is stabilised if there exist positive constants $\tau, \mu$ and positive definite matrix $P \in \Re^{(n+n_c)\times(n+n_c)}$ such that the following matrix inequality holds

$$\begin{bmatrix} \hat{\Lambda}_{11} & \Pi_M & \Pi_N \\ \star & -\mu I & \mu I \\ \star & \star & -Sym(\mathcal{J}) - \mu I \end{bmatrix} < 0 \tag{17}$$

where

$$\begin{aligned} &\hat{\Lambda}_{11} = \begin{bmatrix} P\mathcal{A}_\psi + \mathcal{A}_\psi^T P + \tau\tilde{\xi} I & P\begin{bmatrix} I \\ 0 \end{bmatrix} \\ \star & -\tau I \end{bmatrix} \\ &\Pi_M = \begin{bmatrix} \widetilde{\mathcal{M}} \\ 0 \end{bmatrix}, \Pi_N = \begin{bmatrix} P^T\widetilde{\mathcal{N}}^T \\ 0 \end{bmatrix} \end{aligned} \tag{18}$$

**Proof.** Considering the closed-loop system (15), for any $\mathcal{X}_1(t) = [x_1^T(t) \quad x_{c1}^T(t)]^T$ and $\mathcal{X}_2(t) = [x_2^T(t) \quad x_{c2}^T(t)]^T$ we have

$$\begin{aligned} &\left\| \mathcal{A}_{cl,\Delta}\mathcal{X}_1(t) + \begin{bmatrix} \phi(x_1(t), u(t)) \\ 0 \end{bmatrix} - \mathcal{A}_{cl,\Delta}\mathcal{X}_2(t) - \begin{bmatrix} \phi(x_2(t), u(t)) \\ 0 \end{bmatrix} \right\|_2 \le \|\mathcal{A}_{cl,\Delta}\|_2 \|\mathcal{X}_1(t) - \mathcal{X}_2(t)\|_2 + \\ &\|\phi(x_1(t), u(t)) - \phi(x_2(t), u(t))\|_2 \end{aligned} \tag{19}$$

since $\phi(x, u)$ satisfies the Assumption 3, Lipschitz condition implies that

$$\begin{aligned} &\|\phi(x_1(t), u(t)) - \phi(x_2(t), u(t))\|_2 \\ &\quad < \xi\|[x_1(t) \quad x_2(t)]\|_2 \end{aligned} \tag{20}$$

it can be easily obtained that

$$\begin{aligned} &\xi\|[x_1(t) \quad x_2(t)]\|_2 \le \xi \left\| \begin{bmatrix} x_1(t) & x_2(t) \\ x_{c1}(t) & x_{c2}(t) \end{bmatrix} \right\|_2 = \\ &\xi\|\mathcal{X}_1(t) - \mathcal{X}_2(t)\|_2. \end{aligned} \tag{21}$$

Since matrix $\mathcal{A}_{cl,\Delta}$ have bounded elements, there exist constants $M_1 > 0$ and $M_2 > 0$ such that $\mathcal{A}_{cl,\Delta} \le M_1$. Substituting (21) into inequality (19), it implies that

$$\begin{aligned} &\left\| \mathcal{A}_{cl,\Delta}\mathcal{X}_1(t) + \begin{bmatrix} \phi(x_1(t), u(t)) \\ 0 \end{bmatrix} - \mathcal{A}_{cl,\Delta}\mathcal{X}_2(t) - \right. \\ &\left. \begin{bmatrix} \phi(x_2(t), u(t)) \\ 0 \end{bmatrix} \right\|_2 \le (M_1 + \xi)\|\mathcal{X}_1(t) - \mathcal{X}_2(t)\|_2 \end{aligned} \tag{22}$$

this yields that $\Phi(\mathcal{X}(t), t)$ is Lipschitz in $\mathcal{X}(t)$.

Define $\Phi_{\mathcal{X}}(\mathcal{X}(t), t) = \mathcal{A}_{cl,\Delta}\mathcal{X}(t) + \begin{bmatrix} \phi(x(t), u(t)) \\ 0 \end{bmatrix}$ a continuous function mapping from a set $\Re_\epsilon = \{(t, \mathcal{X}): 0 \le t \le a \ and \ \|\mathcal{X} - \mathcal{X}_0\| \le b\}$ to $\Re^{n+n_c}$. $\Phi(\mathcal{X}(t), t)$ is bounded on $\Re_\epsilon$ with upper bound $M_2 > 0$. It follows from Lemma 1 and Lemma 2 that, the solution of (15) is given by

$$\mathcal{X}(t) = \mathcal{X}_*(t^q/\Gamma(q+1)) \tag{23}$$

where $\mathcal{X}_*(v)$ satisfies the following differential equation

$$\begin{aligned} &d\mathcal{X}_*(v)/dv = \mathcal{A}_{cl,\Delta}\mathcal{X}_*(v) + \Xi(\mathcal{X}_*(v), u_*(v)), \\ &\mathcal{X}_*(0) = [x_0^T \quad x_{c0}^T]^T \end{aligned} \tag{24}$$

with

$$\begin{aligned} &\mathcal{X}_*(v) = \mathcal{X}\left(t - (t^q - v\Gamma(q+1))^{1/q}\right), \\ &x_*(v) = x\left(t - (t^q - v\Gamma(q+1))^{1/q}\right), \\ &x_{c*}(v) = x_c\left(t - (t^q - v\Gamma(q+1))^{1/q}\right), \\ &u_*(v) = u\left(t - (t^q - v\Gamma(q+1))^{1/q}\right), \\ &\Xi(\mathcal{X}_*(t), u_*(t)) = \begin{bmatrix} \phi(x_*(t), u_*(t)) \\ 0 \end{bmatrix} \end{aligned} \tag{25}$$

Consider a candidate Lyapunov function for (24) as follows

$$V(v) = \mathcal{X}_*^T(v) P \mathcal{X}_*(v) \tag{26}$$

where $P$ is a symmetric positive definite matrix. Time derivative of candidate function is calculated as

$$\begin{aligned} \frac{dV(v)}{dv} &= \dot{\mathcal{X}}_*^T(v) P \mathcal{X}_*(v) + \mathcal{X}_*^T(v) P \dot{\mathcal{X}}_*(v) \\ &= \left(\mathcal{A}_{cl,\Delta}\mathcal{X}_*(v) \right. \\ &\quad \left. + \Xi(\mathcal{X}_*(v), u_*(v))\right)^T P \mathcal{X}_*(v) \\ &\quad + \mathcal{X}_*^T(v) P \left(\mathcal{A}_{cl,\Delta}\mathcal{X}_*(v) \right. \\ &\quad \left. + \Xi(\mathcal{X}_*(v), u_*(v))\right) \end{aligned} \tag{27}$$

Introducing the vector $Z = [\mathcal{X}_*^T(v) \ z_1^T]^T$, where $z_1 = \phi(x_*(t), u_*(t))$, and with this in mind that $\mathcal{A}_{cl,\Delta}$ is similar to $\mathcal{A}_{cl,\Delta}^T$ the equation (27) can be rewritten as

$$\begin{aligned} &\frac{dV(v)}{dv} \\ &= Z^T \begin{bmatrix} \mathcal{A}_\psi P + P\mathcal{A}_\psi^T + sym\{\widetilde{\mathcal{M}}\Delta(\sigma)\widetilde{\mathcal{N}}P\} & P\begin{bmatrix} I_n \\ 0 \end{bmatrix} \\ \star & 0 \end{bmatrix} Z \end{aligned} \tag{28}$$

According to direct Lyapunov approach, the stability conditions for the system (24) is $V(v) > 0$ and $dV(v)/dv < 0$. Equation (26) shows that $V(v)$ is positive, and the second condition holds if $dV(v)/dv$ defined in (28) be negative.

It follows from (11) that

$$z_2^T z_2 < \xi^2 \mathcal{X}_*^T(v)\mathcal{X}_*(v) \tag{29}$$

Rearranging inequalities (29) with respect to $Z$, yields

$$Z^T \begin{bmatrix} \xi^2 I & 0 \\ \star & -I \end{bmatrix} Z \ge 0 \tag{30}$$

Applying S-Procedure on $dV(v)/dv < 0$ defined in (28), and (30), and also defining $\tilde{\xi} = \xi^2$, it can be obtained that

$$Z^T \Sigma Z < 0 \tag{31}$$

where

$$\boldsymbol{\Sigma} = \begin{bmatrix} \mathcal{A}_{\psi}\boldsymbol{P} + \boldsymbol{P}\mathcal{A}_{\psi}^{T} + Sym\{\widetilde{\mathcal{M}}\boldsymbol{\Delta}(\sigma)\widetilde{\mathcal{N}}\boldsymbol{P}\} + \tau\tilde{\xi}\boldsymbol{I} & \boldsymbol{P}\begin{bmatrix}\boldsymbol{I}\\ \boldsymbol{0}\end{bmatrix} \\ \star & -\tau\boldsymbol{I} \end{bmatrix}. \tag{32}$$

Let

$$\begin{aligned} &\mathcal{W} = Sym(\mathcal{J}),\\ &\mathcal{Q} = \mathcal{W}^{-1/2}(\widetilde{\mathcal{M}}^{T} + \widetilde{\mathcal{N}}\boldsymbol{P}) - \mathcal{W}^{\frac{1}{2}}\Delta^{T}(\sigma)\widetilde{\mathcal{M}}^{T}, \end{aligned} \tag{33}$$

we have

$$\begin{aligned} -\mathcal{Q}^{T}\mathcal{Q} \le \boldsymbol{0} \;\Leftrightarrow\; &-\Big(\mathcal{W}^{-\frac{1}{2}}(\widetilde{\mathcal{M}}^{T} + \widetilde{\mathcal{N}}\boldsymbol{P}) \\ &- \mathcal{W}^{\frac{1}{2}}\Delta^{T}(\sigma)\widetilde{\mathcal{M}}^{T}\Big)^{T}\Big(\mathcal{W}^{-\frac{1}{2}}(\widetilde{\mathcal{M}}^{T} \\ &+ \widetilde{\mathcal{N}}\boldsymbol{P}) - \mathcal{W}^{\frac{1}{2}}\Delta^{T}(\sigma)\widetilde{\mathcal{M}}^{T}\Big) < \boldsymbol{0} \\ &\Leftrightarrow -sym\{\widetilde{\mathcal{M}}^{T}\mathcal{W}^{-1}\widetilde{\mathcal{N}}\boldsymbol{P}\} \\ &- \widetilde{\mathcal{M}}\mathcal{W}^{-1}\widetilde{\mathcal{M}}^{T} - \boldsymbol{P}^{T}\widetilde{\mathcal{N}}^{T}\mathcal{W}^{-1}\widetilde{\mathcal{N}}\boldsymbol{P} \\ &+ Sym\{\widetilde{\mathcal{M}}\Delta(\sigma)\widetilde{\mathcal{N}}\boldsymbol{P}\} \\ &+ \widetilde{\mathcal{M}}\big(Sym\{\Delta(\sigma)\} \\ &- \Delta(\sigma)\mathcal{W}\Delta^{T}(\sigma)\big)\widetilde{\mathcal{M}}^{T} \le \boldsymbol{0} \\ &\Rightarrow -Sym\{\widetilde{\mathcal{M}}\mathcal{W}^{-1}\widetilde{\mathcal{N}}\boldsymbol{P}\} \\ &- \widetilde{\mathcal{M}}\mathcal{W}^{-1}\widetilde{\mathcal{M}}^{T} - \boldsymbol{P}^{T}\widetilde{\mathcal{N}}^{T}\mathcal{W}^{-1}\widetilde{\mathcal{N}}\boldsymbol{P} \\ &+ Sym\{\widetilde{\mathcal{M}}\Delta(\sigma)\widetilde{\mathcal{N}}\boldsymbol{P}\} \le \boldsymbol{0} \end{aligned} \tag{34}$$

it follows from Lemma 4 that $Sym\{\boldsymbol{\Delta}(\sigma)\} - \boldsymbol{\Delta}^{T}(\sigma)\boldsymbol{\mathcal{W}}\boldsymbol{\Delta}(\sigma) > 0$, and the following inequality holds

$$\begin{aligned} &Sym\{\widetilde{\mathcal{M}}\Delta(\sigma)\widetilde{\mathcal{N}}\boldsymbol{P}\} \le Sym\{\widetilde{\mathcal{M}}\mathcal{W}^{-1}\widetilde{\mathcal{N}}\boldsymbol{P}\} + \\ &\boldsymbol{P}^{T}\widetilde{\mathcal{N}}^{T}\mathcal{W}^{-1}\widetilde{\mathcal{N}}\boldsymbol{P} + \widetilde{\mathcal{M}}\mathcal{W}^{-1}\widetilde{\mathcal{M}}^{T}. \end{aligned} \tag{35}$$

Inequality (35) is equivalent to that there exist $\mu > 0$ such that

$$\begin{aligned} &Sym\{\widetilde{\mathcal{M}}\Delta(\sigma)\widetilde{\mathcal{N}}\boldsymbol{P}\} \le Sym\{\widetilde{\mathcal{M}}\mathcal{W}^{-1}\widetilde{\mathcal{N}}\boldsymbol{P}\} + \\ &\boldsymbol{P}^{T}\widetilde{\mathcal{N}}^{T}\mathcal{W}^{-1}\widetilde{\mathcal{N}}\boldsymbol{P} + \widetilde{\mathcal{M}}(\mathcal{W}^{-1} + \mu^{-1})\widetilde{\mathcal{M}}^{T}. \end{aligned} \tag{36}$$

which is equivalent to that there exist $\mu > 0$ such that

$$\begin{aligned} &sym\{\widetilde{\mathcal{M}}\Delta(\sigma)\widetilde{\mathcal{N}}\boldsymbol{P}\} \le \\ &[\widetilde{\mathcal{M}} \quad \boldsymbol{P}^{T}\widetilde{\mathcal{N}}^{T}]\begin{bmatrix} \mathcal{W}^{-1} + \mu^{-1}\boldsymbol{I} & \mathcal{W}^{-1} \\ \mathcal{W}^{-1} & \mathcal{W}^{-1} \end{bmatrix}\begin{bmatrix} \widetilde{\mathcal{M}}^{T} \\ \widetilde{\mathcal{N}}\boldsymbol{P} \end{bmatrix} = \\ &[\widetilde{\mathcal{M}} \quad \boldsymbol{P}\widetilde{\mathcal{N}}^{T}]\begin{bmatrix} \mu\boldsymbol{I} & -\mu\boldsymbol{I} \\ -\mu\boldsymbol{I} & \mathcal{W} + \mu\boldsymbol{I} \end{bmatrix}^{-1}\begin{bmatrix} \widetilde{\mathcal{M}}^{T} \\ \widetilde{\mathcal{N}}\boldsymbol{P} \end{bmatrix}. \end{aligned} \tag{37}$$

Considering the equation (32) we rewrite the equation (37) as follows

$$\begin{aligned} &\begin{bmatrix} Sym\{\widetilde{\mathcal{M}}\boldsymbol{\Delta}(\sigma)\widetilde{\mathcal{N}}\boldsymbol{P}\} & \boldsymbol{0} \\ \star & \boldsymbol{0} \end{bmatrix} \le \\ &\begin{bmatrix} \widetilde{\mathcal{M}} & \boldsymbol{P}^{T}\widetilde{\mathcal{N}}^{T} \\ \boldsymbol{0} & \boldsymbol{0} \end{bmatrix}\begin{bmatrix} \mu\boldsymbol{I} & -\mu\boldsymbol{I} \\ -\mu\boldsymbol{I} & \mathcal{W} + \mu\boldsymbol{I} \end{bmatrix}^{-1}\begin{bmatrix} \widetilde{\mathcal{M}}^{T} & \boldsymbol{0} \\ \widetilde{\mathcal{N}}\boldsymbol{P} & \boldsymbol{0} \end{bmatrix}. \end{aligned} \tag{38}$$

Substituting (38) into inequality (31), and applying Schur complement completes the proof. ■

Since $\boldsymbol{\mathcal{A}}_{\psi}$ and $\widetilde{\boldsymbol{\mathcal{N}}}$ containing varying terms, is multiplied by $\boldsymbol{P}$ the inequality (17) is bilinear matrix inequality (BMI). To deal with this issue, the following theorem investigates the consensus problem of system (1) in term of LMI (linear matrix inequality).

**Theorem 2.** The output feedback controller (14) solves the stability problem of the system (1) with $0 < \alpha < 1$, if there exist positive constants $\tau$, $\mu$ and positive definite matrices $\boldsymbol{P}_{u} \in \Re^{n\times n}$, $\boldsymbol{P}_{d} \in \Re^{n_c\times n_c}$ and matrices $\mathfrak{A} \in \Re^{n_c\times n_c}, \mathfrak{B} \in \Re^{n_c\times p}, \mathfrak{C} \in \Re^{m\times n_c}, \mathfrak{D} \in \Re^{m\times p}$ such that the following matrix inequality holds

$$\begin{bmatrix} \Lambda_{11} & \Pi_{M} & \Pi_{N} \\ \star & -\mu\boldsymbol{I} & \mu\boldsymbol{I} \\ \star & \star & -Sym(\mathcal{J}) - \mu\boldsymbol{I} \end{bmatrix} < \boldsymbol{0} \tag{39}$$

where

$$\begin{aligned} &\Lambda_{11} = \begin{bmatrix} \Upsilon_{11} + \tau\tilde{\xi}\boldsymbol{I} & \boldsymbol{P}\begin{bmatrix}\boldsymbol{I}\\ \boldsymbol{0}\end{bmatrix} \\ \star & -\tau\boldsymbol{I} \end{bmatrix} \\ &\Upsilon_{11} = \begin{bmatrix} \lambda_{11} & \lambda_{12} \\ \star & \lambda_{22} \end{bmatrix},\ \Pi_{M} = \begin{bmatrix} \widetilde{\mathcal{M}} \\ \boldsymbol{0} \end{bmatrix},\ \Pi_{N} = [\boldsymbol{q}_1 \quad \boldsymbol{q}_2 \quad \boldsymbol{0}]^{T}, \\ &\boldsymbol{P} = \begin{bmatrix} \boldsymbol{P}_{u} & \boldsymbol{0} \\ \boldsymbol{0} & \boldsymbol{P}_{d} \end{bmatrix} \\ &\lambda_{11} = \mathcal{A}\boldsymbol{P}_{u} + \boldsymbol{P}_{u}\mathcal{A}^{T} + \mathcal{B}\mathfrak{D} + \mathfrak{D}^{T}\mathcal{B}^{T}, \\ &\lambda_{12} = \mathcal{B}\mathfrak{C} + \mathfrak{B}^{T},\ \lambda_{22} = \mathfrak{A} + \mathfrak{A}^{T} \\ &\boldsymbol{q}_1 = \mathcal{N}_1\boldsymbol{P}_{u} + \mathcal{N}_2\mathfrak{D},\ \ \boldsymbol{q}_2 = \mathcal{N}_2\mathfrak{C}, \end{aligned} \tag{40}$$

moreover the controller matrices $\boldsymbol{\mathcal{A}}_c$, $\boldsymbol{\mathcal{B}}_c$, $\boldsymbol{\mathcal{C}}_c$ and $\boldsymbol{\mathcal{D}}_c$ are as follows

$$\begin{aligned} &\mathcal{A}_c = \mathfrak{A}\boldsymbol{P}_{d}^{-1}, \quad \mathcal{B}_c = \mathfrak{B}\boldsymbol{P}_{u}^{-1}\mathcal{C}^{\dagger}, \\ &\mathcal{C}_c = \mathfrak{C}\boldsymbol{P}_{d}^{-1}, \quad \mathcal{D}_c = \mathfrak{D}\boldsymbol{P}_{u}^{-1}\mathcal{C}^{\dagger}. \end{aligned} \tag{41}$$

**Proof.** According to the proof of Theorem 1, the output feedback controller (14) solves the consensus problem of the system (1) if the inequality (17) holds. To deal with multiplication of variables, according to $\boldsymbol{P} = diag(\boldsymbol{P}_{u}, \boldsymbol{P}_{d})$ we expand the matrix $\boldsymbol{P}\boldsymbol{\mathcal{A}}_{cl,\Delta}^{T} + \boldsymbol{\mathcal{A}}_{cl,\Delta}\boldsymbol{P}$

$$\begin{aligned} &\boldsymbol{P}\mathcal{A}_{\psi}^{T} + \mathcal{A}_{\psi}\boldsymbol{P} = \begin{bmatrix} \lambda_{11} & \lambda_{12} \\ \lambda_{21} & \lambda_{22} \end{bmatrix},\ \boldsymbol{P}^{T}\widetilde{\mathcal{N}}^{T} = [\boldsymbol{q}_1 \quad \boldsymbol{q}_2]^{T} \\ &\lambda_{11} = \mathcal{A}\boldsymbol{P}_{u} + \boldsymbol{P}_{u}\mathcal{A}^{T} + \mathcal{B}\mathcal{D}_c\mathcal{C}\boldsymbol{P}_{u} + \boldsymbol{P}_{u}\mathcal{C}^{T}\mathcal{D}_c^{T}\mathcal{B}^{T}, \\ &\lambda_{12} = \mathcal{B}\mathcal{C}_c\boldsymbol{P}_{d} + \boldsymbol{P}_{u}\mathcal{C}^{T}\mathcal{B}_c^{T}, \lambda_{21} = \mathcal{B}_c\mathcal{C}\boldsymbol{P}_{u} + \boldsymbol{P}_{d}\mathcal{C}_c^{T}\mathcal{B}^{T} \\ &\lambda_{22} = \mathcal{A}_c\boldsymbol{P}_{d} + \boldsymbol{P}_{d}\mathcal{A}_c^{T} \\ &\boldsymbol{q}_1 = \mathcal{N}_1\boldsymbol{P}_{u} + \mathcal{N}_2\mathcal{D}_c\mathcal{C}\boldsymbol{P}_{u},\ \ \boldsymbol{q}_2 = \mathcal{N}_2\mathcal{C}_c\boldsymbol{P}_{d}. \end{aligned} \tag{42}$$

Now, with the change of the variables as

$$\begin{aligned} &\mathfrak{A} = \mathcal{A}_c\boldsymbol{P}_{d}, \quad \mathfrak{B} = \mathcal{B}_c\mathcal{C}\boldsymbol{P}_{u}, \\ &\mathfrak{C} = \mathcal{C}_c\boldsymbol{P}_{d}, \quad \mathfrak{D} = \mathcal{D}_c\mathcal{C}\boldsymbol{P}_{u}, \end{aligned} \tag{43}$$

equation (42) can be rewritten as

$$\begin{aligned} &\lambda_{11} = \mathcal{A}\boldsymbol{P}_{u} + \boldsymbol{P}_{u}\mathcal{A}^{T} + \mathcal{B}\mathfrak{D} + \mathfrak{D}^{T}\mathcal{B}^{T}, \\ &\lambda_{12} = \mathcal{B}\mathfrak{C} + \mathfrak{B}^{T}, \lambda_{21} = \mathfrak{B} + \mathfrak{C}^{T}\mathcal{B}^{T} \\ &\lambda_{22} = \mathfrak{A} + \mathfrak{A}^{T} \\ &\boldsymbol{q}_1 = \mathcal{N}_1\boldsymbol{P}_{u} + \mathcal{N}_2\mathfrak{D},\ \ \boldsymbol{q}_2 = \mathcal{N}_2\mathfrak{C}, \end{aligned} \tag{44}$$

which completes the proof. ■

**Corollary 1.** [26] Consider fractional-order system (1) without nonlinear term. The output dynamic controller makes the system in (1) asymptotically stable if there exist positive definite Hermitian matrix $\boldsymbol{P} = \boldsymbol{P}^{*}$ in the form of

$$\boldsymbol{P} = diag(\boldsymbol{P}_{u}, \boldsymbol{P}_{d})$$

with $\boldsymbol{P}_{u} \in \mathfrak{C}^{n\times n}$, $\boldsymbol{P}_{d} \in \mathfrak{C}^{n_c\times n_c}$, and matrices $\boldsymbol{\mathcal{T}}_i$, $i = 1, \dots, 4$ and a real constant $\mu > 0$ such that the following LMI constraint become feasible:

$$\begin{bmatrix} \tilde{\Lambda}_{11} & \widetilde{\boldsymbol{M}} & \Lambda_{13} \\ \star & -\mu\boldsymbol{I} & \mu\boldsymbol{I} \\ \star & \star & -Sym(\mathcal{J}) - \mu\boldsymbol{I} \end{bmatrix} < \boldsymbol{0}. \tag{45}$$

where

$$\begin{aligned} &\tilde{\Lambda}_{11} = \begin{bmatrix} \tilde{\lambda}_{11} & \tilde{\lambda}_{12} \\ \star & \tilde{\lambda}_{22} \end{bmatrix}, \Lambda_{13} = \begin{bmatrix} \boldsymbol{q}_1 \\ \boldsymbol{q}_2 \end{bmatrix} \\ &\tilde{\lambda}_{11} = \mathcal{A}(r\boldsymbol{P}_{u} + \bar{r}\overline{\boldsymbol{P}_{u}}) + (r\boldsymbol{P}_{d} + \bar{r}\overline{\boldsymbol{P}_{d}})^{T}\mathcal{A}^{T} + \mathcal{B}\boldsymbol{T}_4 + \\ &\boldsymbol{T}_4^{T}\mathcal{B}^{T}, \\ &\tilde{\lambda}_{12} = \mathcal{B}\boldsymbol{T}_3 + \boldsymbol{T}_2^{T}, \tilde{\lambda}_{21} = \boldsymbol{T}_2 + \boldsymbol{T}_3^{T}\mathcal{B}, \tilde{\lambda}_{22} = \boldsymbol{T}_1 + \boldsymbol{T}_1^{T} \\ &\boldsymbol{q}_1 = (r\boldsymbol{P}_{u} + \bar{r}\overline{\boldsymbol{P}_{u}})^{T}\mathcal{N}_1^{T} + \boldsymbol{T}_4^{T}\mathcal{N}_2^{T}, \boldsymbol{q}_2 = \boldsymbol{T}_3^{T}\mathcal{N}_2^{T} \\ &\boldsymbol{\theta} = (1 - q)\pi/2. \end{aligned} \tag{46}$$

The controller matrices $\boldsymbol{\mathcal{A}}_c$, $\boldsymbol{\mathcal{B}}_c$, $\boldsymbol{\mathcal{C}}_c$ and $\boldsymbol{\mathcal{D}}_c$ can be obtain as follows

$$\begin{aligned} &\mathcal{A}_c = \mathfrak{A}\boldsymbol{P}_{d}^{-1}, \\ &\mathcal{B}_c = diag(\mathfrak{b}_1\bar{p}_u^{-1}\widetilde{\mathcal{C}}^{\dagger}, \dots, \mathfrak{b}_N\bar{p}_u^{-1}\widetilde{\mathcal{C}}^{\dagger}), \\ &\mathcal{C}_c = \mathfrak{C}\boldsymbol{P}_{d}^{-1}, \end{aligned} \tag{47}$$

$\mathcal{D}_c = \boldsymbol{diag}(\mathfrak{d}_1 \bar{p}_u^{-1} \tilde{\mathcal{C}}^{\dagger}, \dots, \mathfrak{d}_N \bar{p}_u^{-1} \tilde{\mathcal{C}}^{\dagger})$.

**Corollary 2.** Proposed methods in Theorem 1, Theorem 2 and Corollary 1 are applicable to the certain form of FO-LTI system (1) by solving the inequalities $\widehat{\boldsymbol{\Lambda}}_{11} < 0$, $\boldsymbol{\Lambda}_{11} < 0$ and $\widehat{\boldsymbol{\Lambda}}_{11}$ respectively.

**Proof.** Assuming $\mathcal{M} = 0$ in the proof Theorem 1, Theorem 2 and Corollary 1, it can be easily obtained.

**Remark 1.** Special case of static output feedback controller for the stabilisation of system (1) can be obtained by solving proposed LMIs for $n_c = 0$.

## IV. Simulation

### A. Example 1

We consider the following non-linear fractional-order system which is available in [10]

$$\begin{aligned} \boldsymbol{D}^{q}\boldsymbol{x}(\boldsymbol{t}) &= \tilde{\mathcal{A}}\boldsymbol{x}(\boldsymbol{t}) + \tilde{\mathcal{B}}\boldsymbol{u}(\boldsymbol{t}) + \boldsymbol{\phi}(\boldsymbol{x}(\boldsymbol{t}), \boldsymbol{u}(\boldsymbol{t})) \\ \boldsymbol{y}(\boldsymbol{t}) &= \mathcal{C}\boldsymbol{x}(\boldsymbol{t}) \end{aligned} \tag{48}$$

with the fractional-order $q = 0.9$ and

$$\begin{aligned} &\mathcal{A} = \begin{bmatrix} \mathbf{0} & \mathbf{1} \\ \mathbf{2} & -\mathbf{6} \end{bmatrix}, \ \mathcal{B} = [\mathbf{1} \quad \mathbf{0.5}]^{T}, \mathcal{C} = [\mathbf{1} \quad \mathbf{1}], \\ &\boldsymbol{\phi}(\boldsymbol{x}(\boldsymbol{t}), \boldsymbol{u}(\boldsymbol{t})) = \begin{bmatrix} \boldsymbol{sin}(\boldsymbol{x}_2) \\ -\boldsymbol{sin}(\boldsymbol{x}_1) + \mathbf{0.5}\, \boldsymbol{sin}(\boldsymbol{x}_2 \boldsymbol{u}(\boldsymbol{t})) \end{bmatrix} \end{aligned} \tag{49}$$

and the uncertainty parameters with

$$\boldsymbol{M} = \begin{bmatrix} \mathbf{0.5} & \mathbf{1} \\ -\mathbf{0.4} & \mathbf{0.2} \end{bmatrix}, \boldsymbol{N}_1 = \begin{bmatrix} \mathbf{0.5} & \mathbf{1.5} \\ \mathbf{0} & \mathbf{0.5} \end{bmatrix}, \boldsymbol{N}_2 = \begin{bmatrix} \mathbf{1} \\ -\mathbf{0.5} \end{bmatrix}^{T}, \boldsymbol{J} = \boldsymbol{I}_2 \tag{50}$$

The time response of the system (48) without control input and $x_0 = [-0.3 \quad 0.3]^T$, is demonstrated in Fig. 1 which shows that states are not convergent and the system (48) is not asymptotically stable.

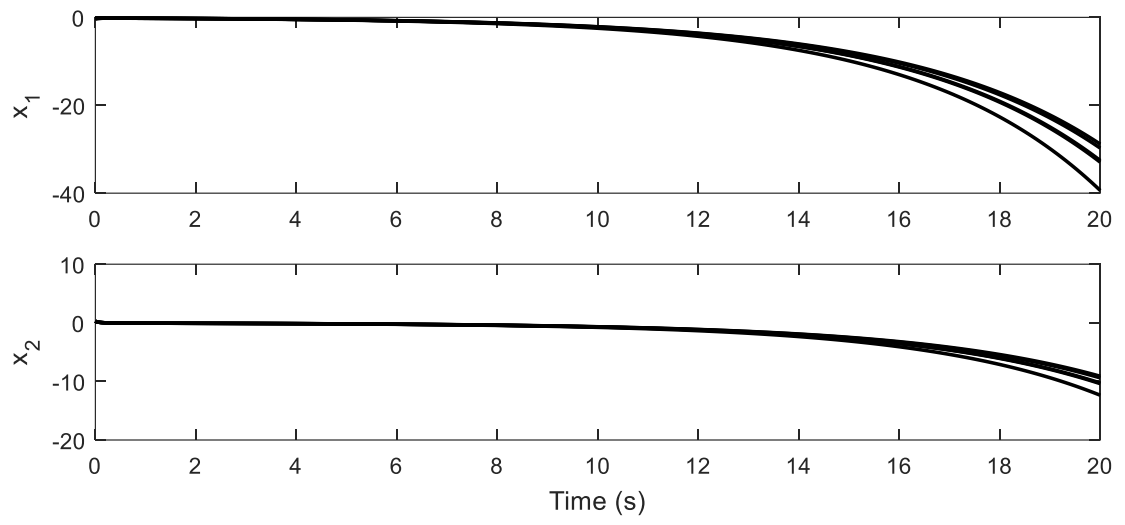

*Fig. 1 Time response of 5 random systems (48) with $u(t) = 0$*

Using Theorem 2, dynamic out feedback controllers of arbitrary orders parameters that stabilise the unstable nonlinear system (48), tabulated in Table 1. Time responses of uncertain nonlinear system (48) via controllers resulted in Table 1 are illustrated in Fig. 2. Results show that unstable system is stabilisable, even with lower orders of dynamic output feedback controllers and all the states asymptotically converge to zero. Nevertheless, comparing the dynamic controller result with the static one indicates that oscillation and settling time of the response of the nonlinear system via dynamic feedback controllers are better than the static output feedback.

*Table 1 controller parameters obtained by Theorem 2 for the system (48)*

| $n_c$ | $A_c$ | $B_c$ | $C_c$ | $D_c$ |
|---|---|---|---|---|
| 0 | 0 | 0 | 0 | $-1.6$ |
| 1 | $-1.3$ | $-2.8$ | 0.6 | $-2.3$ |
| 2 | $\begin{bmatrix} -2.3 & 0 \\ 0.3 & -1.2 \end{bmatrix}$ | $\begin{bmatrix} -0.1 \\ -1.6 \end{bmatrix}$ | $\begin{bmatrix} 0.2 \\ -0.6 \end{bmatrix}$ | $-2$ |

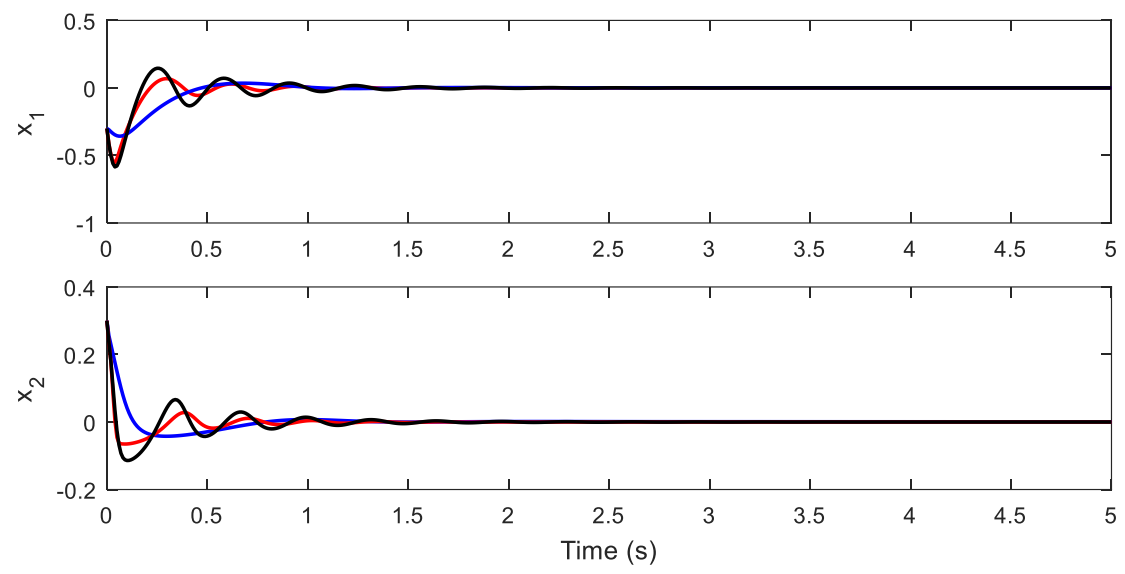

*Fig. 2 Time response of closed-loop system defined in (48) via obtained controllers in Table 1, with $n_c = 2$ (blue), $n_c = 1$ (red), and ordinary static output feedback controller (black).*

To study the robustness of the proposed method in Theorem 2, the resulted dynamic output feedback controller resulted through the Theorem 2 with $n_c = 2$ is utilized for 50 random systems defined in (48). The output time responses depicted in Fig. 3. Results ensure that the proposed method for the control of uncertain systems is reliable and control protocol is effectively robust for the positive definite uncertainty defined in *(7)* to *(9)*.

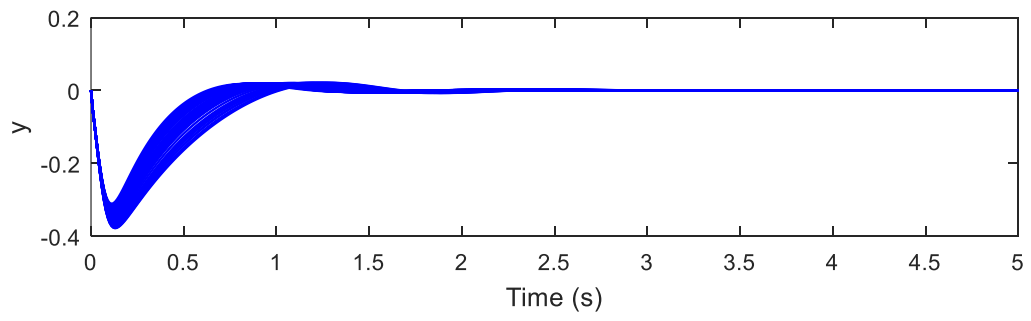

*Fig. 3 Output time response of 50 random system defined in (48) via dynamic controller proposed in Theorem 2 with $n_c = 2$.*

### B. Example 2

The system (48) with the following system matrices is considered

$$\begin{aligned} &\mathcal{A} = \begin{bmatrix} \mathbf{0} & \mathbf{1} \\ \mathbf{2} & -\mathbf{6} \end{bmatrix}, \mathcal{B} = [\mathbf{1} \quad \mathbf{0.5}]^{T}, \mathcal{C} = \begin{bmatrix} \mathbf{1} & \mathbf{2} \\ \mathbf{0.5} & \mathbf{1} \end{bmatrix}, \\ &\boldsymbol{\phi}(\boldsymbol{x}(\boldsymbol{t}), \boldsymbol{u}(\boldsymbol{t})) = \begin{bmatrix} \boldsymbol{sin}(\boldsymbol{x}_2) \\ -\boldsymbol{sin}(\boldsymbol{x}_1) + \mathbf{0.5}\, \boldsymbol{sin}(\boldsymbol{x}_2 \boldsymbol{u}(\boldsymbol{t})) \end{bmatrix} \end{aligned} \tag{51}$$

and the uncertainty parameters are

$$\boldsymbol{M} = \begin{bmatrix} \mathbf{0.5} & \mathbf{1} \\ -\mathbf{0.4} & \mathbf{0.2} \end{bmatrix}, \boldsymbol{N}_1 = \begin{bmatrix} \mathbf{0.5} & \mathbf{1.5} \\ \mathbf{0} & \mathbf{0.5} \end{bmatrix}, \boldsymbol{N}_2 = \begin{bmatrix} \mathbf{1} & \mathbf{1} \\ -\mathbf{0.5} & \mathbf{1.5} \end{bmatrix}^{T}, \boldsymbol{J} = \boldsymbol{I}_2 \tag{52}$$

*Table 2 controllers obtained by Theorem 2 for the system of Example 2*

| $n_c$ | $A_c$ | $B_c$ | $C_c$ | $D_c$ |
|---|---|---|---|---|
| 0 | 0 | 0 | 0 | $\begin{bmatrix} -0.4 \\ 0.1 \end{bmatrix}$ |
| 1 | $-1.4$ | 0.1 | $\begin{bmatrix} -0.3 \\ 0 \end{bmatrix}$ | $\begin{bmatrix} -1.6 \\ 0.3 \end{bmatrix}$ |
| 2 | $\begin{bmatrix} -1.3 & 0 \\ 0 & -1.3 \end{bmatrix}$ | $\begin{bmatrix} 0.4 \\ 0.3 \end{bmatrix}$ | $\begin{bmatrix} 0.1 & 0 \\ -0.5 & -0.3 \end{bmatrix}$ | $\begin{bmatrix} -1.5 \\ 0.2 \end{bmatrix}$ |

In this example we condidered a system with rank defficiency. It's noteworthy that for the abovementioned unstable system no observer-based feedback controller, using the method proposed in [10] can be designed because of the rank constraint on the system matrix $C$. According to the Theorem 2 it can be conculuded that the system (51) is stabilisable by the dynamic output feedback controllers with arbitrary orders tabulated in Table 2.

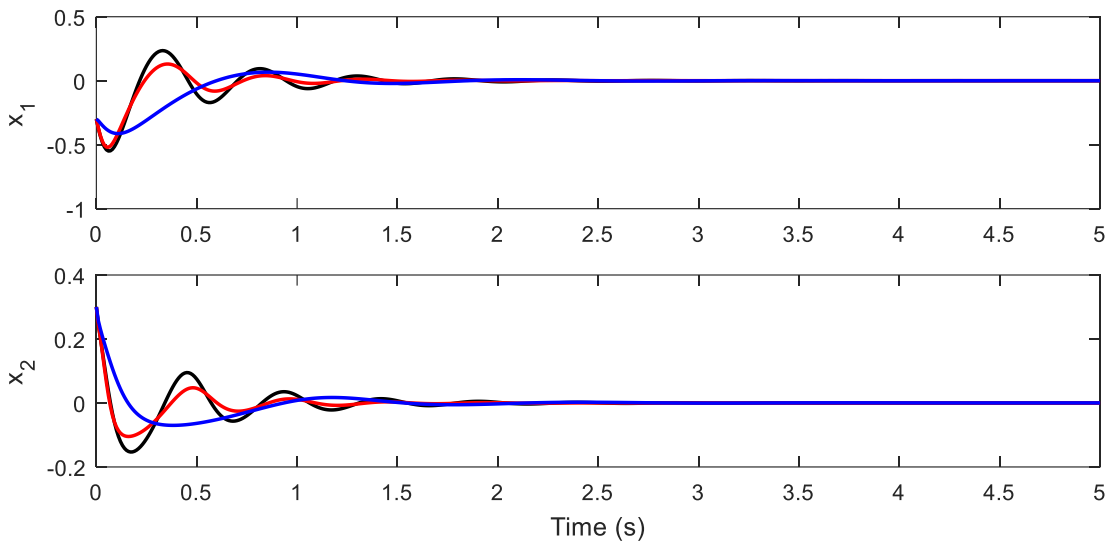

*Fig. 4 Time response of closed-loop system defined in (51) via obtained controllers in Table 2, with $n_c = 2$ (blue), $n_c = 1$ (red), and ordinary static output feedback controller (black).*

Results show that all the states asymptotically converge to zero. Also the effectiveness of proposed dynamic output feedback controller is clear in higher orders which responses have slighter oscillation and shorter settling time. The robustness of the proposed dynamic output feedback controller is studied for the 50 random system (51) via controller in the Table 2 with $n_c = 2$, where the results are depicted in Fig. 5.

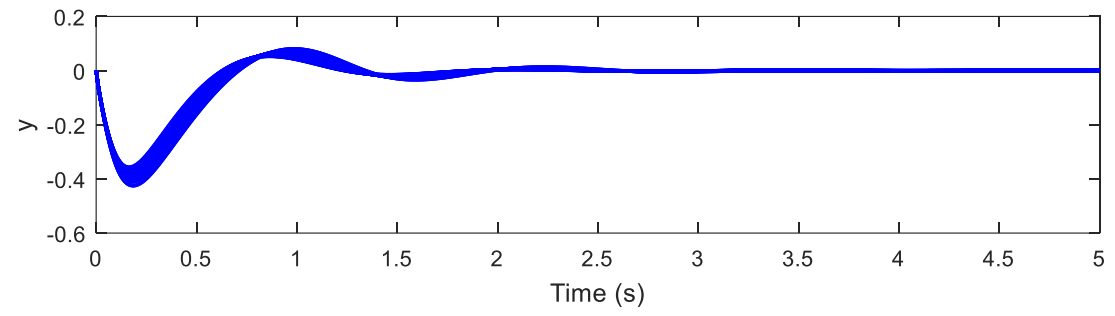

*Fig. 5 Output time response of 50 random system with defined in (51) via dynamic controller proposed in Theorem 2 with $n_c = 2$.*

## V. Conclusion

In this paper, first fixed-order dynamic output feedback controller has been applied to a class of uncertain fractional-order nonlinear system. Then the sufficient conditions for the robust stability of the nonlinear fractional-order system with dynamic output feedback controller with predetermined order, through the direct Lyapunov approach, are derived. The dynamic output feedback controller benefits are accessible, which its order can be set as ideal value in order to reach the desire performance. Note that there are no limitative constraints on the state space matrices of system and the most complete form of dynamic output feedback controller strategy is assumed in our design procedure. Moreover, the result of robust stabilisation is presented in term of LMI, which is straightforward to be utilised. Eventually, some numerical examples are presented to illustrate the effectiveness and advantages of the proposed method.